%% file: egen5.tex
\input  amstex
\input amsppt.sty
\magnification1200
\vsize=23.5truecm
\hsize=16.5truecm
\NoBlackBoxes

\input gmacro2.tex

\def\ul{\underline}

\document

\medskip

\topmatter
\title
Spectral results for mixed problems and fractional 
 elliptic operators 
\endtitle
\author Gerd Grubb \endauthor
\affil
{Department of Mathematical Sciences, Copenhagen University,
Universitetsparken 5, DK-2100 Copenhagen, Denmark.}
E-mail {\tt grubb\@math.ku.dk}
\endaffil
\abstract
 One purpose of the paper is to show Weyl type spectral
asymptotic formulas for pseudodifferential operators $P_a$ of order $2a$, with type and
factorization index $a\in{\Bbb R}_+$ when restricted to a compact set
 with smooth  boundary. The $P_a$ include fractional powers of the Laplace
operator and of variable-coefficient strongly elliptic 
differential operators. Also the regularity of
eigenfunctions is described.

The other purpose is to improve the knowledge of realizations $A_{\chi
,\Sigma _+}$ in $L_2(\Omega )$ of mixed problems for
second-order strongly elliptic symmetric differential operators $A$ on
a bounded smooth set $\Omega \subset {\Bbb R}^n$. Here the boundary
$\partial\Omega =\Sigma $ is partioned smoothly into $\Sigma =\Sigma
_-\cup \Sigma _+$, the Dirichlet
condition $\gamma _0u=0$ is imposed on $\Sigma _-$,  and
a Neumann or Robin condition $\chi u=0$ is imposed on 
$\Sigma _+$. It is shown that the
Dirichlet-to-Neumann operator $P_{\gamma ,\chi }$ is principally of type $\frac12$ with
factorization index $\frac12$, relative to $\Sigma _+$. The above theory allows a detailed
description of $D(A_{\chi ,\Sigma _+})$ with singular elements
outside of $\ol H^{\frac32}(\Omega )$, and leads to a spectral
asymptotic formula for the Krein resolvent difference $A_{\chi ,\Sigma _+}^{-1}-A_\gamma ^{-1}$.  
\endabstract
\rightheadtext{Spectral results}
\endtopmatter

\head Introduction \endhead

This paper has two parts. After a chapter with preliminaries, we
establish in the first part (Chapter 2) spectral asymptotic
formulas of Weyl type for general Dirichlet realizations of
pseudodifferential operators ($\psi $do's) of type $a>0$, as defined in
Grubb \cite{G13,14}, and
discuss the regularity of eigenfunctions.

In the second part (Chapter 3) we consider mixed boundary value problems
for second-order symmetric strongly elliptic differential operators, 
characterize the domain, and find the spectral
asymptotics of the Krein term (the  difference of the resolvent from the
Dirichlet resolvent) in general variable-coefficient situations,
extending the result of \cite{G11a} for the principally Laplacian
case. This includes showing that the relevant Dirichlet-to-Neumann operator
fits into the calculus of the first part.

On Chapter 2: A typical example of the $\psi $do's $P_a$ of type $a>0
$ and order $2a$ that  we treat is the $a$'th power of
the Laplacian $(-\Delta )^a$ on ${\Bbb R}^n$, which is currently of great
interest in probability and finance, mathematical physics and
geometry. Also powers of variable coefficient-operators and much more
general $\psi $dos are included. For the Dirichlet realization
$P_{a,\operatorname{Dir}}$ on a bounded open set $\Omega \subset {\Bbb R}^n$, spectral
studies have mainly been aimed at the fractional Laplacian $(-\Delta
)^a$. In the case of $(-\Delta )^a$,  a Weyl asymptotic formula  was shown already by Blumenthal and
Getoor in \cite{BG59}; recently a refined
asymptotic formula was shown by Frank and Geisinger \cite{FG11}, and
Geisinger gave an
extension to certain other constant-coefficient operators \cite{Ge14}.
The exact domain $D(P_{a,\operatorname{Dir}})$  has not been well described for $a\ge \frac12$, except in integer
cases where the operator belongs to the calculus of Boutet de Monvel
\cite{B71}. Based on a recently published systematic theory \cite{G13} of $\psi $do's of
type $\mu \in {\Bbb C}$ (where those in the Boutet de Monvel calculus
are of type 0), it is now possible to describe domains and
parametrices of operators $D(P_{a,\operatorname{Dir}})$ in an exact
way, when $\Omega $ is smooth. We analyse the sequence of eigenvalues $\lambda _j$
(singular values $s_j$ when the operator is not selfadjoint),
showing that a Weyl asymptotic formula holds in general: 
$$
s_j(P_{a,\operatorname{Dir}})\sim C(P_{a},\Omega )j^{2a/n} \text{ for }j\to\infty ;\tag0.1
$$
moreover we show that the possible eigenfunctions are  in
 $d^aC^{2a}(\comega)$ (in  $d^aC^{2a-\varepsilon }(\comega)$ if $2a\in{\Bbb N}$), where $d(x)\sim
 \operatorname{dist}(x,\partial\Omega )$. The results are generalized to
 operators $P$ of order $m=a+b$ with type and factorization index $a$ ($a,b\in\rp$).

On Chapter 3: The detailed knowledge of $\psi $do's of type $a$ has an application to the classical 
``mixed'' boundary
value problems for a second-order strongly elliptic symmetric  differential
operator $A$ on a smooth bounded set $\Omega \subset {\Bbb R}^n$. Here the boundary condition jumps from a
Dirichlet to a Neumann (or
Robin) condition at the interface of a smooth partition $\Sigma =\Sigma
_-\cup \Sigma _+$ of the boundary $\Sigma =\partial\Omega $; it is also
called the Zaremba problem when $A$ is the Laplacian. 
The 
$L_2$-realization $A_{\chi ,\Sigma _+}$ it defines is less regular than standard
realizations such as the Dirichlet realization $A_\gamma $, but the domain has just been somewhat abstractly
described; it is contained in $\ol H^{\frac32-\varepsilon }(\Omega )$
only (observed by Shamir \cite{S68}),
whereas $D(A_\gamma )\subset \ol H^2(\Omega )$. The resolvent difference
$M=A_{\chi
,\Sigma _+}^{-1}-A_\gamma ^{-1}$ was shown by Birman \cite{B62} to
have eigenvalues satisfying $\mu _j(M)=O(j^{-2/(n-1)})$. The present author studied
$A_{\chi ,\Sigma _+}$ from the point of view of extension theory for elliptic
operators in \cite{G11a} (to which we refer for  more
references to the literature); here we obtained the asymptotic
estimate
$$
\mu _j(M )\sim c(M)j^{-2/(n-1)} \text{ for }j\to\infty ,\tag0.2
$$
in the case where $A$ is principally Laplacian. This was drawing on
the theories of Vishik and Eskin
\cite{E81} and Birman and Solomyak \cite{BS77}, and
other pseudodifferential methods.

We now show that the Dirichlet-to-Neumann
operator $P_{\gamma ,\chi }$ of order 1 on $\Sigma $ associated with $A$ is principally of type
$\frac12$ with factorization index $\frac12$ relative to $\Sigma _+$. In the formulas
connected with the mixed problem, $P_{\gamma ,\chi }$ enters as
truncated to $\Sigma _+$. 
Therefore we can now use the detailed information  on
type $\frac12$ $\psi $do's
to describe the domain of $A_{\chi ,\Sigma _+}$
more precisely, showing how functions $\notin
\ol H^\frac32(\Omega )$ occur. Moreover, using Chapter 2 we can
extend the spectral asymptotic formula (0.2) to the general case
where $A$ has variable coefficients.

\head 1. Preliminaries
\endhead

The notations of \cite{G13,G14b} will be used;
we shall just give a brief summary here.

We consider a Riemannian $n$-dimensional $C^\infty $ manifold $\Omega
_1$ (it can be ${\Bbb R}^n$) and an embedded smooth $n$-dimensional
manifold $\comega$ with boundary $\partial\Omega $ and interior
$\Omega $. For $\Omega _1={\Bbb R}^n$,  $\Omega $ can be $\rnpm=\{x\in
{\Bbb R}^n\mid x_n\gtrless 0\}$; here $(x_1,\dots, x_{n-1})=x'$. In
the general manifold case, $\comega$ is taken compact. For
$\xi \in{\Bbb R}^n$, we
denote $(1+|\xi |^2)^\frac12 =\ang \xi $.
Restriction from $\R^n$ to $\rnp$ resp.\ $\rnm$ (or from
$\Omega _1$ to $\Omega $ resp.\ $\complement\comega$) is denoted $r^+$
resp.\ $r^-$,
 extension by zero from $\rnpm$ to $\R^n$ (or from $\Omega $ resp.\
 $\complement\comega$ to $\Omega _1$) is denoted $e^\pm$. In Chapter
 3, the notation is used for a smooth subset $\Sigma _+$ of an
 $(n-1)$-dimensional manifold $\Sigma $.

A pseudodifferential operator ($\psi $do) $P$ on ${\Bbb R}^n$ is
defined from a symbol $p(x,\xi )$ on ${\Bbb R}^n\times{\Bbb R}^n$ by 
$$
Pu=p(x,D)u=\operatorname{OP}(p(x,\xi ))u 
=(2\pi )^{-n}\int e^{ix\cdot\xi
}p(x,\xi )\hat u\, d\xi =\Cal F^{-1}_{\xi \to x}(p(x,\xi )\hat u(\xi ));\tag1.1
$$  
here $\Cal F$ is the Fourier transform $(\F u)(\xi )=\hat u(\xi
)=\int_{{\Bbb R}^n}e^{-ix\cdot \xi }u(x)\, dx$. The symbol $p$ is
assumed to be such that $\partial_x^\beta \partial_\xi ^\alpha p(x,\xi
)$ is $O(\ang\xi ^{r-|\alpha |})$ for all $\alpha ,\beta $, for some
$r\in{\Bbb R}$ (defining
the symbol class $S^r_{1,0}({\Bbb R}^n\times{\Bbb R}^n)$); then it has
order $r$. The definition of $P$
is carried over to manifolds by use of local coordinates; there are
many textbooks (e.g.\ \cite{G09}) describing this and other rules for
operations with $P$, e.g.\ composition rules.
When $P$ is a $\psi $do on ${\Bbb R}^n$ or $\Omega _1$, $P_+=r^+Pe^+$
denotes its truncation to $\rnp$ resp.\ $\Omega $.

Let $1<p<\infty $  (with $1/p'=1-1/p$), then we define for $s\in{\Bbb
R}$ the Bessel-potential spaces
$$
\aligned
H^s_p(\R^n)&=\{u\in \SD'({\Bbb R}^n)\mid \F^{-1}(\ang{\xi }^s\hat u)\in
L_p(\R^n)\},\\
\dot H^{s}_p(\crnp)&=\{u\in H^{s}_p({\Bbb R}^n)\mid \supp u\subset
\crnp \},\\
\ol H^{s}_p(\rnp)&=\{u\in \D'(\rnp)\mid u=r^+U \text{ for some }U\in
H^{s}_p(\R^n)\};
\endaligned \tag1.2
$$
here $\operatorname{supp}u$ denotes the support of $u$. For $\comega$
compact $\subset\Omega _1$, the definition extends to define $\dot
H^s_p(\comega)$ and $\ol H^s_p(\Omega )$ by use of a finite system of
local coordinates. 
When $p=2$, we get the standard $L_2$-Sobolev spaces, here the lower
index 2 is usually omitted.
 (These and other spaces are thoroughly described in Triebel's book
 \cite{T95}. He  writes $\widetilde H$ instead of $\dot
H$; the present notation stems from H\"ormander's works.)
We also need the
H\"older spaces $C^t$ for $t\in \rp\setminus {\Bbb N}$; when $t\in
 {\Bbb N}_0$, $C^t$ stands for functions with continuous derivatives
 up to order $t$. $\dot C^t(\comega)$ denotes the $C^t$-functions on
 $\Omega _1$ supported in $\comega$. Occasionally, we shall also
 formulate results in
 the H\"older-Zygmund spaces $C^t_*$ for $t\ge 0$, that allow some
 statements to be valid
 for all $t$; they equal $C^t$ when $t\notin {\Bbb N}_0$ and contain
 $C^t$ in the integer cases (more details in \cite{G14b}).
The conventions  $\bigcup_{\varepsilon >0}
H_p^{s+\varepsilon } =
H_p^{s+ 0}$, $\bigcap_{\varepsilon >0} H_p^{s-\varepsilon }=
H_p^{s- 0}$, 
defined in a similar way for the other scales of
spaces, will sometimes be used.

A $\psi $do $P$ is called classical (or polyhomogeneous) when the
symbol $p$ has an asymptotic expansion $p(x,\xi )\sim \sum_{j\in{\Bbb
N}_0}p_j(x,\xi )$ with $p_j$ homogeneous in $\xi $ of degree $m-j$ for
all $j$. Then $P$ has order $m$. One can even allow $m$ to be complex;
then $p\in S^{\operatorname{Re} m}_{1,0}({\Bbb R}^n\times{\Bbb R}^n)$;
the operator and symbol are still said to be of order $m$. 

Here there is an additional definition: $P$ satisfies {\it the $\mu
$-transmission condition} (in short: {\it is of type} $\mu $) for some $\mu
\in{\Bbb C}$  when, in local coordinates,
$$
\partial_x^\beta \partial_\xi ^\alpha p_j(x,-N)=e^{\pi i(m-2\mu -j-|\alpha | )
}\partial_x^\beta \partial_\xi ^\alpha p_j(x,N),\tag1.3
$$
for all $x\in\partial\Omega $, all $j,\alpha ,\beta $, where 
$N$ denotes the interior normal to $\partial\Omega $ at $x$.
The implications of the $\mu $-transmission property 
were a main subject of \cite{G13,G14b}; the mapping properties for such
operators in $C^\infty $-based spaces were shown in H\"ormander
\cite{H85}, Sect.\ 18.2.
  
A special role in the theory is played by the {\it order-reducing
operators}. There is a simple definition of operators $\Xi _\pm^\mu $ on
${\Bbb R}^n$
$$ 
\Xi _\pm^\mu =\operatorname{OP}((\ang{\xi '}\pm i\xi _n)^\mu ) ;
$$
they preserve support
in $\crnpm$, respectively. Here the functions
$(\ang{\xi '}\pm i\xi _n)^\mu $ do not satisfy all the estimates
required for the class $S^{\operatorname{Re}\mu }({\Bbb
R}^n\times{\Bbb R}^n)$, but the operators are useful for some
purposes. There is a more refined choice $\Lambda _\pm^\mu $
that does
satisfy all the estimates, and there is a definition $\Lambda
_\pm^{(\mu )}$ in the manifold situation. These operators 
 define
homeomorphisms for all $s\in{\Bbb R}$ such as 
$$
\aligned
\Lambda^{(\mu )}_+&\colon \dot H^s_p(\comega )\simto
\dot H^{s-\operatorname{Re}\mu }_p(\comega ),\\
\Lambda ^{(\mu )}_{-,+}&\colon \ol H^s_p(\Omega )\simto
\ol H^{s-\operatorname{Re}\mu }_p(\Omega );
\endaligned \tag1.4
$$
 here $\Lambda ^{(\mu )}_{-,+}$ is short for $r^+\Lambda ^{(\mu
 )}_{-}e^+$, suitably extended to large negative $s$ (cf.\ Rem.\ 1.1
 and Th.\ 1.3 in \cite{G13}).

The following special spaces introduced by H\"ormander are particularly adapted to $\mu
$-transmission operators $P$: 
$$
\aligned
H^{\mu (s)}_p(\crnp)&=\Xi _+^{-\mu }e^+\ol H_p^{s-\operatorname{Re}\mu
}(\rnp),\quad  s>\operatorname{Re}\mu -1/p',\\
H^{\mu (s)}_p(\comega)&=\Lambda  _+^{(-\mu )}e^+\ol H_p^{s-\operatorname{Re}\mu
}(\Omega ),\quad  s>\operatorname{Re}\mu -1/p',\\
\Cal E_\mu (\comega)&=e^+\{u(x)=d(x)^\mu v(x)\mid v\in C^\infty
(\comega)\};
\endaligned\tag1.5
$$
namely, $r^+P$ (of order $m$) maps them into $\ol H_p^{s-\operatorname{Re}m}(\rnp)$, $\ol
 H_p^{s-\operatorname{Re}m}(\Omega)$ resp.\ $C^\infty (\comega)$ (cf.\
 \cite{G13} Sections 1.3, 2, 4), and they appear as domains of
 elliptic realizations of $P$. In the third line, $\operatorname{Re}\mu
 >-1$ (for other $\mu $, cf.\ \cite{G13}) and $d(x)$ is a $C^\infty
 $-function positive on $\Omega $ and vanishing to
 order 1 at $\partial\Omega $, e.g.\
 $d(x)=\operatorname{dist}(x,\partial\Omega )$ near $\partial\Omega $.
 One has that $H^{\mu
 (s)}_p(\comega)\supset \dot H_p^s(\comega)$, and that the distributions are
 locally in $H^s_p$   on $\Omega $, but at the boundary they in general have a
 singular behavior. More details are given in \cite{G13, G14b}.

\head 2. Spectral results for Dirichlet realizations of type $a$ operators
\endhead

\subhead 2.1 Dirichlet realizations of type $a$ operators\endsubhead 

Consider a Riemannian $n$-dimensional $C^\infty $ manifold $\Omega
_1$ ($n\ge 2$)  
and an embedded compact $n$-dimensional
$C^\infty $-manifold $\comega$ with boundary $\partial\Omega $ and interior
$\Omega $. We consider an elliptic  pseudodifferential operator on
$\Omega _1$ with the following properties  explained in detail in 
\cite{G13}:

\proclaim{Assumption 2.1}
Let $a\in\rp$. $P_a$ is  a classical elliptic $\psi $do on $\Omega _1$
of order $2a$, which relative to $\Omega $ satisfies the
$a$-transmission condition and has the factorization index $a$. 
\endproclaim

For example, $P_a$ can be the $a$-th power of a strongly
elliptic second-order differential operator on $\Omega _1$, in
particular $(-\Delta )^a$, or it can be the $a/m$-th power of a
properly elliptic differential operator of even order $2m$, but also
other operators are allowed.   (We call such operators ``fractional
elliptic'', because they share important properties with 
the fractional Laplacian.) 

As in \cite{G13}, we define the Dirichlet realization
$P_{a,\operatorname{Dir}}$ in $L_2(\Omega )$ as the operator acting like $r^+P_a$ with domain
$$
D(P_{a,\operatorname{Dir}})=\{u\in \dot H^a(\comega)\mid r^+P_au\in
L_2(\Omega )\}.\tag 2.1
$$
Then according to \cite{G13}, Sections 4-5,
$$
D(P_{a,\operatorname{Dir}})=H^{a(2a)}(\comega)=\Lambda _+^{(-a)}e^+\ol H^a(\Omega ).\tag2.2
$$

We recall from \cite{G13}:

\proclaim{Lemma 2.2} For $1<p<\infty $, $s>a-1/p'$,  
the spaces
$H_p^{a(s)}(\comega)$ satisfy:
$$
H_p^{a(s)}(\comega)=\Lambda _+^{(-a)}e^+\ol H^{s-a}_p(\Omega )
\cases =\dot
H_p^{s}(\comega),\text{ if }s-a\in \,]-1/p',1/p[ \,,\\
\subset\dot
H_p^{s-0}(\comega),\text{ if }s=a+1/p,\\
\subset d^ae^+\ol H_p^{s-a}(\Omega ) +\dot H_p^{s}(\comega),\text{ if }s-a-1/p
\in\rp\setminus{\Bbb N},\\
\subset d^ae^+\ol H_p^{s-a}(\Omega ) +\dot H_p^{s-0}(\comega),\text{ if
}s-a-1/p\in \Bbb N.\endcases \tag2.3
$$

Moreover, 
$$H_p^{a(s)}(\comega)\subset \dot H_p^a(\comega), \text{ when
}s-a\ge 0.\tag2.4$$
\endproclaim

\demo{Proof} The equalities in (2.3) come from the definition of
$H_p^{a(s)}(\comega)$, and the inclusions are special cases of
\cite{G13} Th.\ 5.4. For the last statement, we note that when $s-a\ge
0$, $e^+\ol H_p^{s-a}(\Omega )\subset e^+L_p(\Omega )$, which is
mapped into $\dot H_p^a(\comega)$ by
$\Lambda _+^{(-a)}$.\qed 
\enddemo

In the case where $P_a$ is strongly elliptic, i.e., the principal symbol
$p_{a,0}(x,\xi )$ satisfies
$$
\operatorname{Re}p_{a,0}(x,\xi )\ge c|\xi |^{2a},
$$
with $c>0$,  we can describe
$D(P_{a,\operatorname{Dir}})$ in a different way: 

Modifying $\Omega _1$ at a distance from $\comega$ if necessary, we
can assune $\Omega _1$ to be compact without boundary. 
Then it is well-known that  $P_a$ satisfies a G\aa{}rding inequality for $u\in C^\infty (\Omega _1)$:
$$
\operatorname{Re}(P_au,u)_{L_2(\Omega _1)} \ge c_0\|u\|^2_{H^{ a}(\Omega _1)}-k\|u\|^2_{L_2(\Omega _1)},\tag2.5
$$
with $c_0>0$, $k\in {\Bbb R}$ (cf.\ e.g.\ \cite{G09}, Ch.\ 7), besides the inequality 
$$
|(P_au,v)_{L_2(\Omega _1)} |\le C \|u\|_{H^{ a}(\Omega _1)}\|v\|_{H^{ a}(\Omega _1)}.
$$
(In the case of $(-\Delta )^a$ on ${\Bbb
R}^n$, $\Omega \subset {\Bbb R}^n$, there is a slightly different
formulation: For general $P_a$ one would here require $x$-estimates of
the symbol to be uniform on the noncompact set ${\Bbb R}^n$; see e.g.\
\cite{G11c} for the appropriate version of the G\aa{}rding inequality. One can also include
this case by replacing ${\Bbb R}^n\setminus \Omega $ by a suitable
compact manifold.)

Define the sesquilinear
form $s_0$ on $C_0^\infty (\Omega )$ by
$$
s_0(u,v)=(r^+P_au,v)_{L_2(\Omega )}=(P_au,v)_{L_2(\Omega _1)}, \text{
for } u,v\in C_0^\infty (\Omega );
$$
 it extends by closure to a bounded sesquilinear form $s(u,v)$ on $\dot
H^a(\comega)$, to which the inequality (2.5) extends.
The Lax-Milgram
construction applied to $s(u,v)$ (cf.\ e.g.\ \cite{G09}, Ch.\ 12) leads to an operator $S$ which acts 
like $r^+P_a\colon \dot H^a(\comega)\to \ol H^{\,-a}(\Omega )$, with
domain consiting of the functions that are 
mapped into $L_2(\Omega )$; this is exactly $P_{a,\operatorname{Dir}}$
as in (2.1), (2.2). Here both $S$ and $S^*$ are lower bounded, with
lower bound $>-k$ (they are in fact sectorial), hence have $\{\lambda
\in{\Bbb C}\mid \lambda \le -k\}$ in their resolvent sets.

When $P_a$ moreover is symmetric, $P_{a,\operatorname{Dir}}$
 is the Friedrichs extension of $(r^+P_a)|_{C_0^\infty (\Omega )}$.

In the case
of $P_a=(-\Delta )^a$, some authors for
precision call this $P_{a,\operatorname{Dir}}$ the ``restricted fractional Laplacian'', see e.g.\ Bonforte,
Sire and Vazquez 
\cite{BSV14}, in order to distinguish it from the ``spectral
fractional Laplacian'' defined as the $a$'th power of the Dirichlet realization
of
$-\Delta $.

\subhead 2.2 Regularity of eigenfunctions \endsubhead 

The possible eigenfunctions have a certain smoothness:

\proclaim{Theorem 2.3} Let $P_a$ satisfy Assumption {\rm 2.1}.

If $0$ is an eigenvalue of $P_{a,\operatorname{Dir}}$, its associated eigenfunctions are in
$\E_a(\comega)$.  

When $a\in \rp\setminus{\Bbb N}$, then the
eigenfunctions $u $ of $P_{a,\operatorname{Dir}} $ associated with
nonzero eigenvalues $\lambda $ lie in $ d^aC^{2a}
(\comega)$ if $2a\notin {\Bbb N}$, in
$ d^aC^{2a-\varepsilon }
(\comega)$ (any $\varepsilon >0$)
if $2a\in{\Bbb N}$; they are also in $C^\infty (\Omega )$.

When $a\in{\Bbb N} $,  the
eigenfunctions $u $ of $P_{a,\operatorname{Dir}} $ associated with an
eigenvalue $\lambda  $ lie in 
$\{u\in 
C^\infty (\comega)\mid \gamma _0u=\gamma _1u=\dots =\gamma
_{a-1}u=0\}$ (equal to $\E_a(\comega)$ in this case).
\endproclaim

\demo{Proof}  (In some of the formulas here, the extension by zero $e^+$
is tacitly understood.) 
When $\lambda  $ is an
eigenvalue, the associated eigenfunctions $u $
are nontrivial solutions of 
$$
r^+P_au =\lambda  u .\tag2.6
$$
If $\lambda  =0$, then $u \in\E_a(\comega)$, since the right-hand
side in (2.6) is in $C^\infty (\comega)$, and we can apply \cite{G13}
Th.\ 4.4.

Now let 
$\lambda  \ne 0$. 
When $a\in{\Bbb N}$, we are in a well-known standard elliptic case (as
treated e.g.\ in \cite{G96}, Sect.\ 1.7); the eigenfunctions are in
$C^\infty (\comega)$ as well as in $\E_a(\comega)$, and 
$\E_a(\comega)$ is the described subset of $ C^\infty (\comega)$.

Next, consider the case $a\in
\rp\setminus {\Bbb N}$.

To
begin with, we know that $u \in \dot H^a(\comega)$ (from (2.1)). We shall use
the well-known general embedding properties for $p,p_1\in \,]1,\infty [\,$:
$$
\dot H_p^{a}(\comega)\subset e^+L_{p_1}(\Omega ),\text{ when }\tfrac
n{p_1}\ge \tfrac n p - a,  \quad \dot H_p^{a}(\comega)\subset
\dot C^0(\comega )\text{ when }a>\tfrac np.\tag2.7 
$$

If $a> \frac n2$, we have already that $\dot H^a(\comega)\subset \dot C^0 (\comega
)$, so (2.6) has right-hand side in $C^0(\comega )$, and we can
go on with solution results in H\"older spaces; 
this will be done further below.

If $a\le\frac n2$, we make a finite number of iterative steps to reach
the information $u \in C^0(\comega )$, as follows: 
Define $p_0,p_1, p_2,\dots$, with $p_0=2$ and $q_j=\frac n {p_j}$ for
all the relevant $j$,
such that
$$
q_j=q_{j-1}-a \text{ for }j\ge 1.
$$
 This means that $q_j=q_0-ja$; we stop the sequence at $j_0$ the first
 time we reach a  $q_{j_0}\le 0$. As a first step, we note that $u \in \dot
 H^a(\comega)\subset e^+L_{p_1}(\Omega )$ implies $u \in
 H_{p_1}^{a(2a)}(\comega)$ by
 \cite{G13} Th.\ 4.4, and then by (2.4),  $u\in \dot
 H^a_{p_1}(\comega)$. 
In the next step we use the embedding $\dot
 H^a_{p_1}(\comega)\subset e^+L_{p_2}(\Omega )$ to conclude in a
 similar way that
 $u \in \dot H^a_{p_2}(\comega)$, and so on. If $q_{j_0}<0$, we have
 that $\frac n{p_{j_0}}<a$, so $u \in \dot
 H^a_{p_{j_0}}(\comega)\subset \dot C^0(\comega)$. If $q_{j_0}=0$, the
 corresponding $p_{j_0}$ would be $+\infty $, and we see at least that
 $u \in e^+L_{p}(\Omega )$ for any large $p$; then one step more
 gives that $u \in \dot C^0(\comega)$.
 
The rest of the argumentation relies on H\"older estimates, as in \cite{G13},
Section 7, or still more efficiently by \cite{G14b}, Section 3. By the
regularity results there,
$$
u \in C^0(\comega) \implies u \in C^{a(2a)}_*(\comega)\subset
e^+d^a C^{a}(\comega)+\dot C^{2a-0 }(\comega)\subset
e^+C^{a}(\comega).$$ 
Next, $u \in C^{a}(\comega)$ implies $$
u \in C^{a(3a)}_*(\comega)\subset
e^+d^a C^{2a}_*(\comega)+\dot C^{3a\,(-\varepsilon )}(\comega)\subset
e^+d^aC^{2a\,(-\varepsilon )}(\comega) 
$$
where $(-\varepsilon )$ is active if $2a\in{\Bbb N}$. Moreover, by the
ellipticity of $P_a-\lambda  $ on $\Omega _1$, $u $ is $C^\infty $ on
the interior  $\Omega $.
\qed

\enddemo

The fact that an eigenfunction in $\dot H^a(\comega)$ is in $L_\infty
(\Omega )$ was shown for $P_a=(-\Delta )^a$ with $0<a<1$ by Servadei and
Valdinoci \cite{SV13} by a completely different method.

\example{Remark 2.4}
For  $P_a=(-\Delta )^a$ it has been shown by
Ros-Oton and Serra (see \cite{RS14}) that an
eigenfunction $u $ cannot have $u /d^a$ vanishing identically on
$\partial\Omega $. This implies that the regularity of $u $ cannot be
improved all the way up to $\E_a(\comega)$, when $\lambda  \ne 0$,
$a\in\rp\setminus{\Bbb N}$. For if $u $ were in $
\E_a(\comega)$, it would also lie in $C^\infty (\comega)$ (since
$r^+P_au =\lambda  u $ would lie there). Now it is easily checked that
$C^\infty (\comega)\cap
\E_a(\comega)=\dot C^\infty (\comega)$ when $a\in{\Bbb R}_+\setminus
{\Bbb N}$, where the functions vanish to order $\infty $ at the
boundary. In particular, $u /d^a$ would be zero on $\partial\Omega $, contradicting
$u \ne 0$. 
\endexample

The theorem extends without difficulty to operators of order $m=a+b$
considered in $H^s_p$-spaces:

\proclaim{Theorem 2.5} Let $P$ be of type $a>0$ with factorization
index $a$, and of order $m=a+b$, $b>0$. Let $1<p<\infty $, and define
$P_{\operatorname{Dir}}$ as the operator from $H_p^{a(m)}(\comega)$ to
$L_p(\Omega )$ acting like $r^+P$. If $0$ is an eigenvalue, the
associated eigenfunctions are in $\E_a(\comega)$. If $\lambda \ne 0$
is an eigenvalue, the associated eigenfunctions are in $d^a
C^m(\comega)$ (in $d^aC^{m-\varepsilon }(\comega)$ if $m$ is integer).
\endproclaim

\demo{Proof} The zero eigenfunctions are solutions with a $C^\infty $
right-hand side, hence lie in $\E_a(\comega)$ by \cite{G13} Th.\ 4.4.

Now let u be an eigenfunction associated with an eigenvalue $\lambda
\ne 0$. In view
of (2.4), we have $u\in \dot H_p^a(\comega)$. Using (2.7), we find by
application of the regularity result of \cite{G13} Th.\ 4.4, by a
finite number of iterative steps as in the proof of Theorem 2.3, that $u\in \dot H_{p_1}^a,\dot
H_{p_2}^a,\dots $ with increasing $p_j$'s, until we reach $u\in C^0(\comega)$. 

Now we can apply the H\"older results from \cite{G13}, \cite{G14b}; this
goes most efficiently by \cite{G14b} Th.\ 3.2 $2^\circ$ and Th.\ 3.3 for
H\"older-Zygmund spaces:
$$
r^+Pu\in \ol C^{t}_*(\Omega) \implies u\in C^{a(m+t)}_*(\comega)\subset
d^ae^+\ol C_*^{m+t-a}(\Omega )+\dot C_*^{m+t\,(-\varepsilon )}(\comega),\tag2.8
$$
$t\ge 0$,  where $(-\varepsilon )$ is active if $m+t-a\in{\Bbb N}$.

If $b> a$, there are two steps:
$$
u\in C^0(\comega) \implies u\in C^{a(a+b)}_*(\comega)\subset
e^+d^a \ol C_*^{b}(\Omega)+\dot C_*^{a+b(-\varepsilon ) }(\comega)\subset
e^+\ol C_*^{a}(\Omega).$$ 
Next, $u\in \ol C_*^{a}(\Omega)$ implies $$
u\in C_*^{a(m+a)}(\comega)\subset
e^+d^a \ol C^{m}_*(\Omega)+\dot C_*^{m+a\,(-\varepsilon )}(\comega)\subset
e^+d^aC^{m\,(-\varepsilon )}(\comega), 
$$
where $(-\varepsilon )$ is active if $m\in{\Bbb N}$. 

If $b\le a$, we need a finite number of steps, such as:
$$
u\in C^0(\comega) \implies u\in C^{a(a+b)}_*(\comega)\subset
e^+d^a \ol C_*^{b}(\Omega)+\dot C_*^{a+b(-\varepsilon ) }(\comega)\subset
e^+\ol C_*^{b}(\Omega),$$ 
where we use that $a+b-\varepsilon >b$ for small $\varepsilon $. Next, $u\in \ol C_*^{b}(\Omega)$ implies $$
u\in C_*^{a(m+b)}(\comega)\subset
e^+d^a \ol C^{2b}_*(\Omega)+\dot C_*^{a+2b\,(-\varepsilon )}(\comega)\subset
e^+\ol C^{\min\{2b,a\}}_*(\Omega), 
$$
where we use that $a+2b-\varepsilon >\min\{2b,a\}$ for small $\varepsilon $. If
$2b\ge a$, we end the proof as above. If not, we estimate again, now
arriving at the exponent $\min\{3b,a\}$, etc., continuing until we
reach $kb\ge a$; then the proof
is completed as above.
\qed

\enddemo

\subhead 2.3 Spectral asymptotics \endsubhead

We shall now study spectral asymptotic estimates for our operators. We first recall some
notation and basic rules.

As in \cite{G84} we denote by ${\frak C}_p(H,H_1)$ the  $p$-th Schatten class 
consisting of the compact operators $B$ from a Hilbert space $H$ to
another $H_1$ such that
$(s_j(B))_{j\in{\Bbb N}}\in \ell_p({\Bbb N})$. Here the $s$-numbers, or
singular values, are defined as $ s_j(B)=\mu
_j(B^*B)^\frac12$, where $\mu _j(B^*B)$ denotes the $j$-th positive
eigenvalue of $B^*B$, arranged nonincreasingly and repeated according
to multiplicities. 
The so-called weak Schatten class
consists of the compact operators $B$ such that
$$
s_j(B)\le Cj^{-1/p}\text{ for all }j;\text{ we set } {\bold
N}_p(B)=\sup_{j\in{\Bbb N}}s_j(B)j^{1/p}.\tag2.9
$$
The notation $\frak S_{(p)}(H,H_1)$ was used in
\cite{G84} for this space; instead we here use the name  $\frak
S_{p,\infty }(H,H_1)$ (as in \cite{G14a} and in other works).
The indication $(H,H_1)$ is replaced by $(H)$ if $H=H_1$; it can be
omitted when it is clear from the context. One has that $\frak
S_{p,\infty }\subset {\frak C}_{p+\varepsilon } $ for any $\varepsilon
>0$. They are linear spaces.

We recall (cf.\ e.g.\ \cite{G84} for details and references) that ${\bold N}_p(B)$ is
a quasinorm on $\frak S_{p,\infty }$, with a good control over the
behavior under summation.
Recall also that
$$
\frak S_{p,\infty }\cdot \frak S_{q,\infty }\subset \frak S_{r,\infty }, \text{ where }r^{-1}=p^{-1}+q^{-1},\tag2.10
$$
and
$$
s_j(B^*)=s_j(B),\quad
s_j(EBF)\le \|E\|s_j(B)\|F\|,\tag2.11
$$
when $E\colon H_1\to H_3$ and $F\colon H_2\to H$ are bounded linear maps between Hilbert
spaces. 

Moreover, we recall that when $\Xi $ is a bounded open subset
of ${\Bbb R}^m$ and reasonably regular, or is a compact smooth
$m$-dimensional manifold with boundary, then the injection $H^t(\Xi
)\hookrightarrow L_2(\Xi )$ is in $\frak S_{m/t,\infty }$ when $t>0$. It
follows that when $B$ is a linear operator in $L_2(\Xi )$ that is
bounded from $L_2(\Xi )$ to $H^t(\Xi )$, then $B\in \frak
S_{m/t,\infty }$, with
$$
{\bold N}_{m/t}(B)\le C\|B\|_{\Cal L(L_2(\Xi ),H^t(\Xi ))}.\tag2.12
$$

Recall also the Weyl-Ky Fan perturbation result: 
$$ 
s_j(B)j^{1/p}\to C_0,\quad s_j(B')j^{1/p}\to 0 \implies
s_j(B+B')j^{1/p}\to C_0, \text{ for }j\to\infty.\tag2.13
$$

We shall moreover use Laptev's result \cite{L81}: When $P$ is a
classical $\psi $do of order $t<0$ on a closed $m$-dimensional
manifold $\Xi _1$ with a smooth subset $\Xi $, $m\ge 2$, then
$$
1_{\Xi _1\setminus\Xi }P 1_{\Xi }\in \frak S_{(m-1)/t,\infty };\tag2.14
$$
in fact it has a Weyl-type asymptotic formula of that order.

Results on the spectral behavior of compositions of $\psi $do's of
negative order interspersed with
functions with jumps were shown in \cite{G11b}, see in particular Th.\
4.3 there. 
We need to supply this
result
 with a statement allowing a zero-order factor of the form
of a sum of a pseudodifferential and a singular Green operator (in the
Boutet de Monvel calculus); as
functions with jumps we here just take $1_\Omega $. 

\proclaim{Theorem 2.6}
Let $M_\Omega $ be an operator on $\comega$ composed of ${l}\ge 1$ factors
$P_{j,+}$ formed of
classical pseudodifferential operators $P_j$ on $\Omega _1$ of negative orders
$-t_j$ and  truncated to $\comega$, $j=1,\dots,l$, and one factor $Q_++G$
(placed somewhere between them), where $Q$ is classical of order $0$
and $G$ is a singular Green operator on $\comega$ of order and class
$0$: 
$$
M_\Omega =P_{1,+}\dots P_{l_0,+}(Q_++G)P_{l_0+1,+}\dots P_{l,+}.\tag2.15
$$
Let $t=t_1+\dots+t_{l}$, and let $m(x,\xi )$ be the product of the
principal $\psi $do symbols on $\Omega _1$:
$$
m(x,\xi )=p_{1,0}(x,\xi )\dots q_0(x,\xi )\dots p_{l,0}(x,\xi ).
$$
Then $M_\Omega $ has the spectral behavior:
$$
s_j(M_\Omega )j^{t/n} \to c(M_\Omega )^{t/n}\text{ for }j\to\infty ,\tag2.16
$$
where 
$$
c(M_\Omega ) =
\tfrac1{n(2\pi )^{n}}\int_{\Omega  }\int_{|\xi
|=1}
\big(m(x,\xi )^*m(x,\xi )
\big)^{n/2t}
\,d\omega (\xi ) dx
\tag2.17
$$
\endproclaim

\demo{Proof} By Th.\ 4.3 of \cite{G11b} with interspersed functions
of the form $1_{\Omega }$, the statement
holds if $Q=1$ and $G=0$, so the new thing is to include nontrivial cases
of $Q$ and $G$. We can assume that $l_0\ge 1$. For the contribution from $Q$ we write
$$
P_{l_0,+}Q_+=r^+P_{l_0}e^+r^+Qe^+=r^+P_{l_0}Qe^+-r^+P_{l_0}e^-r^-Qe^+.\tag2.18
$$
Here $P_{l_0}Q$ is a $\psi $do of order $-t_{l_0}<0$ with principal
symbol $p_{l_0,0}q_0$, and when $r^+P_{l_0}Qe^+$ is taken into the
original expression, we get an operator of the type treated by Th.\
4.3 of \cite{G11b},
$$
P_{1,+}\dots (P_{l_0}Q)_+P_{l_0+1,+}\dots P_{l,+},\tag2.19
$$
for which the
statement (2.16), (2.17) holds.
For the other term in (2.18), we use that $r^+P_{l_0}e^-$ is the type
of operator covered by the theorem of Laptev \cite{L81} (cf.\ (2.14)), belonging to $\frak S_{(n-1)/t_{l_0},\infty }$,
and    $r^-Qe^+$ is bounded in $L_2$, so in view of the rules (2.10)
and (2.11)
for compositions, the full expression with
this term inserted is in $\frak S_{n/(t+\theta  ),\infty }$ for a
certain $\theta  >0$. The spectral asymptotic estimate obtained for
the term (2.19) is preserved
when we add this term of a better weak Schatten class, in view of (2.13).

The contribution from $G$ will likewise be shown to be in a better
weak Schatten class that the main $\psi $do term; this requres a
deeper effort. Actually, the strategy can be copied from some proofs in
\cite{G14a}, as follows: Consider first the composition of $G$ with just one operator:
$$
M=P_+G,
$$
where $P$ is of order $-t<0$. In local coordinates, we can extend  Th.\ 4.1 in \cite{G14a}
to this operator, writing
$$
\psi P_+G \psi _1=\sum_{k\in{\Bbb N}_0}\psi P_+K_k\Phi ^*_k\psi
_1=\sum_{k\in{\Bbb N}_0}\psi P_+\zeta K_k\Phi ^*_k\psi _1+\sum_{k\in {\Bbb N}_0}\psi P_+K_k(1-\zeta )\Phi ^*_k\psi _1,
$$
with Poisson and trace operators $K_k$ and $\Phi ^*_k$ as explained in
 \cite{G14a}, and letting $P_+K_k$ play the role of $K_k$ in the proof
 there. Here $(\psi
 P_+K_k\zeta )^*$ is bounded from $L_2(B_{R,+})$ to $\ol H^t(B'_{R'})$ for a large $R'$, hence lies in $\frak
 S_{(n-1)/t,\infty }$ (by the property of the injection of  $\ol
 H^t(B'_{R'})$ into  $L_2(B'_{R'})$, $B'_{R'}=\{x'\in{\Bbb
 R}^{n-1}\mid |x'|<R'\}$). The proof that the full series $P_+G$ lies in $\frak
 S_{(n-1)/t,\infty }$ goes as in \cite{G14a} (using also that the terms
 with $1-\zeta $ have a smoothing component, and that the series is
 rapidly convergent). Moreover, Cor.\ 4.2 there
 shows how the result is extended to the manifold situation.

When there are several factors in $M$, we need only use that
$P_{j,+}\in \frak S_{n/t_j,\infty }$ for the other factors and apply
the product rule (2.10), and we end with
the information that the full product is in $\frak S_{n/(t+\theta ),\infty }$ for
some $\theta >0$, so that the spectral asymptotics 
remains as that of (2.19), when this term is added on.
\qed
\enddemo

The result extends easily to matrix-formed operators.

Now we can show a spectral asymptotic estimate for $P_{a,\operatorname{Dir}}$. 

\proclaim{Theorem 2.7} Let $P_a$ satisfy Assumption {\rm 2.1}. Assume that $P_{a,\operatorname{Dir}}$ is
invertible, or more generally that $P_{a,\operatorname{Dir}}+c$ is
invertible from $D(P_{a,\operatorname{Dir}})$ to $L_2(\Omega )$ for some
$c\in{\Bbb C}$ (this holds if $P_a$ is strongly elliptic).

The singular values
$s_j(P_{a,\operatorname{Dir}})$ (eigenvalues of
$(P_{a,\operatorname{Dir}}^*P_{a,\operatorname{Dir}})^\frac12$) have the
asymptotic behavior:
$$
s_j(P_{a,\operatorname{Dir}})=C(P_{a,\operatorname{Dir}})j^{2a/n}+o(j^{2a/n}),
\text{ for }j\to\infty ,\tag2.20
$$
where $C(P_{a,\operatorname{Dir}})=C'(P_{a,\operatorname{Dir}})^{-2a/n}$,
defined from the principal symbol $p_{a,0}(x,\xi )$ by: 
$$
C'(P_{a,\operatorname{Dir}})=\tfrac 1{n(2\pi )^n}\int_{\Omega }\int_{|\xi |=1}|p_{a,0}(x,\xi )|^{-n/2a}\,d\omega (\xi )dx.\tag2.21
$$

\endproclaim

\demo{Proof} By Th.\ 4.4 of \cite{G13},  $P_{a,\operatorname{Dir}}$,
acting like $r^+P_a$, has
a parametrix of order $-2a$,$$
R=\Lambda _{+,+}^{(-a)}(\widetilde Q_++G)\Lambda _{-,+}^{(-a)}=r^+\Lambda _+^{(-a)}e^+(r^+\widetilde Qe^++G)r^+\Lambda _-^{(-a)}e^+
; \tag2.22$$
in the last expression, we have written the restriction- and extension-operators out in
detail. In comparison with the formula for $R$ in \cite{G13}, Th.\ 4.4, we have moreover
placed an $r^+$ in front, which is allowed since $R$ maps into a space
of functions supported in $\comega$. (The singular Green operator
component $G$ was missing in some preliminary versions of \cite{G13}.) The
operator is of the form treated in Theorem 2.6, which gives the
asymptotic behavior of the $s$-numbers of $R$:
$$
s_j(R)j^{2a/n}\to c(R)^{2a/n} \text{ for }j\to\infty ;\tag2.23
$$
here $c(R)=C'(P_{a,\operatorname{Dir}})$ defined in (2.21), since the principal symbol of
$\Lambda _+^{(-a)}\widetilde Q\Lambda _-^{(-a)}$ is the inverse of the
principal symbol of $P_a$. 

That $R$ is  parametrix of $r^+P_a=P_{a,\operatorname{Dir}}$ implies that
$$
P_{a,\operatorname{Dir}}R=I-S_1, \text{ where }S_1\colon L_2(\Omega )\to C^\infty (\comega).\tag2.24
$$
 Consider the case where $P_{a,\operatorname{Dir}}$ is invertible; it is
 clearly compact since it maps $L_2(\Omega )$ into $\dot
 H^a(\comega)$. It follows from (2.24) that
$$
P_{a,\operatorname{Dir}}^{-1}=P_{a,\operatorname{Dir}}^{-1}(P_{a,\operatorname{Dir}}R+S_1)=R+S_2,
\quad S_2=P_{a,\operatorname{Dir}}^{-1}S_1, 
$$
where $P_{a,\operatorname{Dir}}^{-1}\in \frak S_{n/a,\infty }$ (since it
maps $L_2(\Omega )$ into $\dot H^a(\comega)$), and
$S_1\in \bigcap_{p>0}\frak S_{p,\infty }$, so $S_2\in \bigcap_{p>0}\frak
S_{p,\infty }$ by (2.10). By (2.13), the spectral asymptotic formula
(2.23) for $R$ will
therefore imply the same spectral asymptotic formula for
$P_{a,\operatorname{Dir}}^{-1}$, so
$$
s_j(P_{a,\operatorname{Dir}}^{-1})j^{2a/n}\to C'(P_{a,\operatorname{Dir}}^{-1})^{2a/n}.
$$
The asymptotic formula can 
 also be writen as the formula (2.20)
for the $s$-numbers of $P_{a,\operatorname{Dir}}$. 

If instead $P_{a,\operatorname{Dir}}+c$ is invertible,
we can write
$$
(P_{a,\operatorname{Dir}}+c)R=I-S_1+cR,
$$
with $S_1$ as in (2.24), and hence
$$
\aligned
(P_{a,\operatorname{Dir}}+c)^{-1}&=(P_{a,\operatorname{Dir}}+c)^{-1}((P_{a,\operatorname{Dir}}+c)R+S_1-cR)\\
&=R+(P_{a,\operatorname{Dir}}+c)^{-1}S_1-c(P_{a,\operatorname{Dir}}+c)^{-1}R.
\endaligned
$$
Here $(P_{a,\operatorname{Dir}}+c)^{-1}S_1\in \bigcap_{p>0}\frak
S_{p,\infty }$ and $c(P_{a,\operatorname{Dir}}+c)^{-1}R\in \frak
S_{n/3a,\infty }$, since   $(P_{a,\operatorname{Dir}}+c)^{-1}\in \frak S_{n/a,\infty }$, and
$R\in \frak S_{n/2a,\infty }$ in view of its spectral behavior shown
above.
Thus $(P_{a,\operatorname{Dir}}+c)^{-1}$ is a perturbation of $R$ by
operators in better weak Schatten classes, and the desired spectral results
follow for $(P_{a,\operatorname{Dir}}+c)^{-1}$ and it inverse
$P_{a,\operatorname{Dir}}+c$.
\qed 
\enddemo

When $P_{a,\operatorname{Dir}}$ is selfadjoint $\ge 0$, its eigenvalue
sequence $\lambda _j, j\in {\Bbb N}$, coincides with
the sequence of $s_j$-values, and
 Theorem 2.7 gives an
asymptotic estimate of the eigenvalues.

In this case, the asymptotic estimate  extends to  arbitrary open sets $\Omega $ (assumed bounded when $\Omega _1={\Bbb
R}^n$), with the
Dirichlet
realization defined by Friedrichs extension of $r^+P_a$ from $C_0^\infty (\Omega )$,
since the eigenvalues can be characterized by the minimax principle,
which gives a monotonicity property in terms of nested open sets.

As mentioned in the introduction, the estimate (2.20) was shown for
the case $P_a=(-\Delta )^a$ by Blumenthal and Getoor in
\cite{BG59}. In this case, a
two-terms asymptotic formula for the $N$'th average of eigenvalues as
$N\to\infty $ was obtained by Frank and Geisinger in \cite{FG11}, and
Geisinger extended the  estimate (2.20) to a larger class of 
constant-coefficient $\psi $do's in \cite{Ge14}.

\example{Remark 2.8} Theorem 2.7 extends straightforwardly to Dirichlet realizations of operators $P$
as in Theorem 2.5; in the proof, the factor $\Lambda ^{(-a)}_{-,+}$ is
replaced by $\Lambda ^{(-b)}_{-,+}$, and $2a$ in the asymptotic
formula is replaced by $m=a+b$.
\endexample

\head 3. Mixed problems for  second-order symmetric strongly elliptic
differential operators \endhead

\subhead 3.1 The Krein resolvent formula \endsubhead

We shall now apply the knowledge of the operators of type $\frac12$ to
the mixed boundary value problem for second-order elliptic
differential operators. The setting is the following:

On a bounded $C^\infty $-smooth open subset $\Omega $ of ${\Bbb R}^n$ with boundary
$\partial\Omega =\Sigma $ we consider a second-order symmetric
differential operator with real coefficients in $C^\infty (\comega)$:
$$
Au= -{\sum}_{j,k=1}^n \partial _{ j} (a_{jk}(x)
  \partial _{ k} u) + a_0(x)u,\tag3.1
$$
here $a_{jk}=a_{kj}$ for all $j,k$. 
$A$ is assumed strongly elliptic, i.e., $\sum_{j,k=1}^n a_{jk}(x)\xi _j\xi _k\ge c_0|\xi
|^2$ for $x\in\comega $, $\xi \in{\Bbb R}^n$, with $c_0>0$.
We denote as usual $u|_{\Sigma }=\gamma _0u$, and consider moreover
the conormal derivative $\nu $ and a Robin variant $\chi $ (both are Neumann-type boundary operators) 
$$
\nu u = {\sum}_{j,k=1}^n n_j \gamma_0 (a_{jk}\partial_{ k} u),\quad \chi u=\nu u-\sigma\gamma _0u;
\tag3.2
$$
 here $\vec
n=(n_1,\dots,n_n)$ denotes the interior unit normal to the boundary, and
 $\sigma$ is a real $C^\infty $-function on $\Sigma $. 
With $\Sigma _+$ denoting a closed $C^\infty $-subset of $\Sigma $, we define
$L_2(\Omega )$-realizations  $A_\gamma $ and $A_{\chi ,\Sigma _+}$ of
 $A$ determined respectively by the boundary conditions:
$$
\aligned
\gamma _0u&=0 \text{ on }\Sigma ,\text{ the Dirichlet condition},
\\
\chi u&=0 \text{ on }\Sigma _+,\;\gamma _0u=0\text{ on
}\Sigma \setminus\Sigma _+,\text{ a  mixed condition}.
\endaligned\tag3.3
$$

It is accounted for in \cite{G11a} that with the domains defined more
precisely by$$\aligned
D(A_{\gamma })&=
\{u\in \ol H^2(\Omega )\mid  \gamma _0u=0\},\\
D(A_{\chi ,\Sigma _+})&=
\{u\in \ol H^1(\Omega )\cap D(\Ama)\mid  \gamma _0u\in \dot H^\frac12( \Sigma _+),\,
\chi u=0\text{ on }\Sigma _+^\circ\},
\endaligned\tag3.4
$$
where $\Ama$ denotes the operator acting like $A$ with domain  $D(\Ama)=\{u\in L_2(\Omega )\mid Au\in L_2(\Omega )\}$, the
operators $A_\gamma $ and $A_{\chi ,\Sigma _+}$ are selfadjoint lower
bounded. We can and shall assume that a sufficiently large constant has been
added to $A$ such that both operators have a positive lower bound.

Let
$$
X=\dot H^{-\frac12}(\Sigma _+);\text{ then }X^*=\ol H^\frac12(\Sigma _+^\circ),\tag3.5
$$
with respect to a duality consistent with the $L_2$-scalar product on
$\Sigma _+$. The injection $\inj_X\colon X\hookrightarrow
H^{-\frac12}(\Sigma )$ can be viewed as an ``extension by zero''
$e^+$ (often tacitly understood), and its adjoint $(\inj_X)^*\colon H^{\frac12}(\Sigma )\to \ol
H^{\frac12}(\Sigma _+^\circ)$ is the restriction $r^+$.

Recalling that $\gamma _0$ defines a homeomorphism from
$Z=\operatorname{ker}(\Ama)=\{u\in L_2(\Omega )\mid Au=0\}$ to $H^{-\frac12}(\Sigma )$
with inverse $K_\gamma $ (a Poisson operator), we define
$$
V=K_\gamma X,\quad \gamma _V\colon V\simto X;\tag3.6
$$
here $V$ is a closed subspace of $Z$ (both closed in the $L_2(\Omega
)$-norm), and $\gamma _V$ denotes the restriction of $\gamma _0$ to
$V$. Note that $\gamma _V^{-1}$ acts like $K_\gamma $ on $X$; it is also
denoted $K_{\gamma ,X}$ in \cite{G11a}. We denote by
${\inj_V}$ the injection of $V$ into $Z$, its adjoint is the
orthogonal projection $\pr_V$ of $Z$ onto $V$. Let us
moreover introduce the relevant Dirichlet-to-Neumann operators
$$
P_{\gamma ,\nu }=\nu K_\gamma ,\quad  P_{\gamma ,\chi }=\chi K_\gamma =P_{\gamma ,\nu }-\sigma ;\tag3.7
$$
they are pseudodifferential operators of order 1 on $\Sigma $, both formally selfadjoint.

The following Krein
resolvent formula was shown in \cite{G11a}, Sect.\ 4.1:

\proclaim{Proposition 3.1} For the realizations of $A$ defined above,
$$
A_{\chi ,\Sigma _+}^{-1}-A_\gamma ^{-1}=\inj_V\gamma _V^{-1}L^{-1}(\gamma _V^{-1})^*\pr_V.\tag3.8
$$
Here $L$ is the (selfadjoint invertible) operator from $X$ to $X^*$ acting like $-r^+P_{\gamma ,\chi }e^+$
and with domain
$$
D(L)=\gamma _0 D(A_{\chi ,\Sigma _+}).
$$
\endproclaim 

It was shown in \cite{G11a} that $D(L)\subset \dot H^{1-\varepsilon
}(\Sigma _+)$ for all $\varepsilon >0$, but that the inclusion does not
hold with $\varepsilon =0$.

Since $L$ acts like $-P_{\gamma ,\chi ,+}$ and is surjective onto $\ol
H^{\frac12}(\Sigma _+^\circ)$, we also have
$$
D(L)=\{\varphi \in \dot H^{1-\varepsilon
}(\Sigma _+) \mid r^+P_{\gamma ,\chi }\varphi \in \ol H^{\frac12}(\Sigma ^\circ_+) \}.\tag3.9
$$
Below we shall improve the knowledge
of the domain by setting $P_{\gamma ,\chi }$ in relation to the types of
operators studied in Chapter 2.

\subhead 3.2 Structure of the Dirichlet-to-Neumann operator
\endsubhead

To study the symbol of $P_{\gamma ,\chi }$ we consider the operators in
a neighborhood $U$ of a point $x_0\in \partial\Omega =\Sigma $, where
local coordinates $x=(x_1\dots,x_n)=(x',x_n)$  are chosen such that $U\cap
\Omega =\{(x',x_n)\mid x'\in B_1, 0<x_n<1\}$ and $U\cap
\partial\Omega =\{(x',x_n)\mid x'\in B_1, x_n=0\}$; $B_1=\{x'\in{\Bbb
R}^{n-1}\mid |\xi '|<1\}$. In these
coordinates, the principal symbol of $A$ at the boundary is a polynomial $$
\aligned
\ul a(x',0,\xi
)&=\sum_{j,k=1}^n\ul a_{jk}(x',0)\xi _j\xi _k=\ul a_{nn}(x',0)\xi
_n^2+2b(x',\xi ')\xi _n+c(x',\xi '),\text{ with}\\
b&=\sum_{j=1}^{n-1}\ul a_{jn}(x')\xi _j,\quad c=\sum_{j,k=1}^{n-1}\ul
a_{jk}(x')\xi _j\xi _k;
\endaligned\tag3.10
$$
the coefficients are real with $\ul a_{jk}=\ul
a_{kj}$. We often write $(x',0)$ as $x'$. Since $A$ is strongly elliptic, $\ul a(x',\xi ',\xi _n)>0$ when
$\xi '\ne 0$,  so the polynomial $\ul a(x',\xi ',\lambda )$ in $\lambda $ has no real roots
when $\xi '\ne 0$. When we set
$$
a'(x',\xi ')=\ul a_{nn}(x')c(x',\xi ')-b(x',\xi ')^2=\sum_{j,k=1}^{n-1}a'_{jk}(x')\xi _j\xi _k,
$$
we therefore have that $a'(x',\xi ')>0$ for $\xi '\in{\Bbb
R}^{n-1}\setminus 0$. The roots of $\ul a(x',\xi ',\lambda )$  equal $\lambda _\pm=\ul a_{nn}^{-1}(-b\pm i\kappa
_0)$, lying respectively in ${\Bbb C}_\pm=\{\lambda \in{\Bbb C}\mid
\operatorname{Im}\lambda \gtrless 0\}$, where we have set
$$
\kappa _0(x',\xi ')=a'(x',\xi ')^{\frac12}>0.\tag3.11
$$
Denote 
$$
\kappa _{\pm}(x',\xi ')=\mp i\lambda _{\pm}=\ul a_{nn}^{-1}(\kappa
_0\pm ib);\tag3.12
$$
then $\ul a$ has the factorization
$$
\ul a(x',\xi ',\xi _n)=\ul a_{nn}(x')(\kappa _+(x',\xi ')+i\xi _n)(\kappa _-(x',\xi ')-i\xi _n),\tag3.13
$$
where $\kappa _+$ and $\kappa _-$ both have positive real part
($=\kappa _0$). This plays a role in standard investigations of
boundary problems. We go on to study the Dirichlet-to-Neumann operators.

The principal symbol-kernel $\tilde k_\gamma (x',x_n,\xi ')$ of $K_\gamma $ is the solution operator
for the semi-homogeneous model problem (with $\varphi $ given in ${\Bbb C}$):
$$
\ul a(x',\xi ',D_n)u(x_n)=0\text{ on }\rp, \quad u(0)=\varphi ;
$$
it is seen from (3.13) that the solution in $L_2(\rp)$ is 
$\varphi e^{-\kappa _+x_n}$, so
$$
\tilde k_\gamma (x',x_n,\xi ')=e^{-\kappa _+x_n}.\tag3.14
$$ 
The conormal derivative for the model problem is
$$
\nu u=\gamma _0(\ul a_{nn}\partial_{x_n}u(x_n)+\sum_{k=1}^{n-1}\ul a_{nk}i\xi _ku(x_n)).
$$
Then the principal
symbol of $P_{\gamma ,\nu  }$ is 
$$
\aligned
p_{\gamma ,\nu  }(x',\xi ')_1&=\gamma _0(\ul
a_{nn}\partial_{x_n}+\sum_{k=1}^{n-1}\ul a_{nk}i\xi _k)e^{\kappa
_+x_n}=
-\ul a_{nn}\kappa _++\sum_{k=1}^{n-1}\ul a_{nk}i\xi _k\\
&=-\ul a_{nn}(-i)\ul a_{nn}^{-1}(-b+i\kappa _0)+ib=-\kappa _0.
\endaligned
$$
Since $P_{\gamma ,\chi }=P_{\gamma ,\nu }-\sigma $ with
$\sigma $ of order 0, $P_{\gamma ,\chi }$ likewise has the principal symbol
$-\kappa _0$.

The important fact that we observe here is that $\kappa _0(x',\xi ')$
is {\it even} in $\xi '$; 
$$
\aligned
\kappa _0(x',-\xi ')&=\kappa _0(x',\xi '), \text{ with
}\\
\partial_{x'}^\beta \partial_{\xi '}^\alpha \kappa _0(x',-\xi
')&=(-1)^{|\alpha |}\partial_{x'}^\beta \partial_{\xi '}^\alpha \kappa
_0(x',\xi ')\text{ for all }\alpha ,\beta ,
\endaligned\tag3.15
$$
(since $c(x',\xi ')$ and $b(x',\xi ')^2$ are clearly even in $\xi '$).
Since $\kappa _0$ is homogeneous of degree 1, it therefore has the
$\frac12$-transmission property with respect to any smooth subset of
$B_1$, satisfying (1.4) with $m=1$, $\mu =\frac12$. 

Moreover, we shall show that it has factorization index $\frac12$ with respect to any
smooth subset of $B_1$: We can take the subset as $B_{1,+}=\{x'\in{\Bbb
R}^{n-1}\mid |x '|<1,x_{n-1}>0 \}$, with $(x_1,\dots,x_{n-2})$ denoted
$x''$. Now we apply the same procedure as above to the polynomial
$a'(x'',0,\xi ')=\kappa _0(x'',0,\xi '',\xi _{n-1})^2$ in $\xi _{n-1}$.
It has a factorization analogously to (3.13):
$$
\kappa _0(x'',0,\xi ')^2=a'_{n-1,n-1}(x'')(\kappa '_+(x'',\xi '')+i\xi _{n-1})(\kappa '_-(x'',\xi '')-i\xi _{n-1}),
$$
where $a'_{n-1,n-1}>0$ and 
$\kappa '_\pm$ have positive real part; here $\kappa '_{\pm}= \mp i
\lambda '_{\pm}$, where $\lambda '_{\pm}$ are the roots of $a'(x'',0,\xi
'',\lambda )$ lying in ${\Bbb C}_\pm$, respectively. It
follows that
$$
\kappa _0(x'',0,\xi ')=a'_{n-1,n-1}(x'')^\frac12(\kappa '_+(x'',\xi '')+i\xi _{n-1})^\frac12(\kappa '_-(x'',\xi '')-i\xi _{n-1})^\frac12,\tag3.16
$$
where $(\kappa '_+(x'',\xi '')+ i\xi _{n-1})^\frac12$ extends
analytically in $\xi _{n-1}$ into ${\Bbb C}_-$ and $(\kappa '_-(x'',\xi '')- i\xi _{n-1})^\frac12$ extends
analytically in $\xi _{n-1}$ into ${\Bbb C}_+$ (in short, are a ``plus-symbol'' resp.\ a
``minus-symbol'', cf.\ \cite{E81}, \cite{G13}). 

Carrying the information
back to $\Omega $ and $\Sigma =\partial\Omega $, we have obtained:

\proclaim{Theorem 3.2} The principal symbol of the
Dirichlet-to-Neumann operator $P_{\gamma ,\chi }$
equals $-\kappa _0(x',\xi ')$ (expressed in local coordinates in
{\rm(3.10)-(3.11)}), negative and elliptic of order $1$. For any smooth subset $\Sigma _+$ of $\Sigma
 $, $\kappa _0$ is of type $\frac12$ and has factorization index $\frac12$
relative to $\Sigma _+$. An explicit factorization in local
coordinates is given in {\rm (3.16)}.
\endproclaim

\subhead 3.3 Precisions on $L$ and $L^{-1}$ 
\endsubhead

Define $L_1$ to be a $\psi $do on $\Sigma $ with symbol $ \kappa
_0(x',\xi ')$, and let $L_0=-P_{\gamma ,\chi }-L_1$. Then since $L $ acts
like $-P_{\gamma ,\chi ,+}$, it acts like $L_{1,+}+L_{0,+}$:
$$
L\varphi = L_{1,+}\varphi +L_{0,+} \varphi , \text{ for }\varphi \in D(L).\tag3.17
$$
Here $L_1$, classical of order 1, is principally equal to $-P_{\gamma ,\chi }$ and $-P_{\gamma ,\nu }$, whereas the operator $L_0$ is a classical $\psi $do of order 0, containing
both the local term $\sigma $ and the nonlocal difference between $P_{\gamma ,\nu
}$ and its principal part. 

As
shown in Theorem 3.2, $L_1$  is of type
$\frac12$ and has factorization index $\frac12$ relative to $\Sigma _+$. Here $L_{1,+}$, when
considered on $\dot H^{1-\varepsilon }(\Sigma _+)$, identifies with the
operator $r^+L_1$ in the homogeneous Dirichlet problem for $L_1$,
going from $\dot H^{1-\varepsilon }(\Sigma _+)$ to $\dot
H^{-\varepsilon }(\Sigma _+)$. It has according to \cite{G13} Th.\ 4.4 a parametrix  $R\colon \ol H^{s-1}(\Sigma
_+^\circ)\to H^{\frac12(s)}(\Sigma _+)$ for $s>\frac12$; here 
$H^{\frac12(s)}(\Sigma _+)=\dot H^s(\Sigma _+)$ for $\frac12<s<1$,
cf.\ (2.3), and $R$ is of the form
$$
R=\Lambda _{+,+}^{(-\frac12)}(\widetilde Q_++G)\Lambda _{-,+}^{(-\frac12)}, \tag3.18
$$
with a $\psi $do $\widetilde Q$ of order and type 0 and a singular
Green operator $G$ of order and class 0.
The parametrix property implies that
$$
\aligned
L_{1,+}R=I-S_1,\quad S_1&\colon \ol H^t(\Sigma _+)\to C^\infty (\Sigma
_+),\text{ for }t>-\tfrac12,\\
RL_{1,+}=I-S_2,\quad S_2&\colon \dot H^{1+t}(\Sigma _+)\to
\E_{\frac12}(\Sigma _+),  \text{ for }-\tfrac12<t<0,\\
\qquad\qquad  S_2&\colon  H^{\frac12(1+t)}(\Sigma _+)\to
\E_{\frac12}(\Sigma _+),\text{ for }t\ge 0.
\endaligned\tag3.19
$$

From (3.17) and the first line in (3.19),
 we have for the difference $S_3$ of $L^{-1}$ and $R$:
$$
S_3=L^{-1}-R=L^{-1}(L_{1,+}R+S_1)-L^{-1}(L_{1,+}+L_{0,+})R=L^{-1}S_1-L^{-1}L_{0,+}R.\tag3.20
$$
 Some properties of $L^{-1}$ can be obtained by considerations similar to
 those in \cite{G11a}:

\proclaim{Proposition 3.3}
The operator $L^{-1}\colon X^*\to X$ extends to an operator $M_0$ that maps continuously
$$
M_0\colon  \ol H^s(\Sigma _+^\circ)\to \dot H^{s+\frac12-\varepsilon }(\Sigma
_+)\text{ for }-1<s\le \tfrac12, \text{ any }\varepsilon >0.
$$
In particular, the closure of $L^{-1}$ in $L_2(\Sigma _+)$ is a continuous operator 
from $L_2(\Sigma _+)$ to $\dot H^{\frac12-\varepsilon }(\Sigma _+)$.

The operators $L^{-1}$ and $M_0$ have the same eigenfunctions (for
nonzero eigenvalues); they
belong to $D(L)$.
\endproclaim  

\demo{Proof} We already know from \cite{G11a} (cf.\ (3.9)) that $L^{-1}$ is
continuous from $X^*=\ol H^\frac12(\Sigma _+^\circ)$ to $\dot H^{1-\varepsilon
}(\Sigma _+)$. Then it has an adjoint $M_0$ (with respect to dualities
consistent with the $L_2(\Sigma _+)$-scalar product) that is
continuous from $\ol H^{\,-1+\varepsilon }(\Sigma _+^\circ)$ to
$\dot H^{-\frac12}(\Sigma _+)$. But since $L^{-1}$ is known to be
selfadjoint (from $X^*$ to $X$, consistently with the $L_2$-scalar product), $M_0$ must be an extension of
$L^{-1}$. Now the asserted continuity for $-1<s\le\frac12$ follows  by
interpolation. For $s=0$ this shows the mapping property of the
$L_2$-closure.

When $\varphi $ is a distribution in $\ol H^{\,-1+\varepsilon }(\Sigma
_+^\circ)$ such that $M_0\varphi =\lambda \varphi $ for some $\lambda \ne 0$,
then since $M_0\varphi \in \ol H^{\,-\frac12+\varepsilon }(\Sigma
_+^\circ)=\dot H^{-\frac12+\varepsilon }(\Sigma _+)$, $\varphi $ lies
there. Next, it follows that $M_0\varphi \in \ol H^{\varepsilon _1}(\Sigma
_+^\circ)=\dot H^{\varepsilon _1}(\Sigma _+)$, and hence $\varphi $ also lies
there. 
Finally, we conclude that $M_0\varphi \in \ol H^{\frac12 +\varepsilon _2}(\Sigma
_+^\circ)$, so that $\varphi $ also lies there. Here $M_0$ coincides
with $L^{-1}$.
\qed
\enddemo

\comment
 $L^{-1}$ extends to an operator $M_0$ that maps continuously 
$
M_0\colon \ol H^{s}(\Sigma _+^\circ)\to \dot H   ^{s+\frac12-\varepsilon }(\Sigma _+)
$ for $-1<s\le \frac12$, any $\varepsilon >0$. In particular, the closure of $L^{-1}$ in
$L_2(\Sigma _+)$ maps $L_2(\Sigma _+)$ into $\dot
 H^{\frac12-\varepsilon }(\Sigma _+)$.
\endproclaim

\demo{Proof} The proof  goes exactly as the proof of
\cite{G11a} Lemma 5.10; it is a simple argumentation using the
selfadjointness of $L^{-1}$, the dual statement of the mapping property
$L^{-1}\colon \ol H^{\frac12}(\Sigma _+^\circ)\to \dot
H^{1-\varepsilon }(\Sigma _+)$, and interpolation.
\qed
\enddemo
\endcomment

We can now find exact information on the domain of $L$:

\proclaim{Theorem 3.4} $L^{-1}$ maps $\ol H^{\frac12}(\Sigma
_+^\circ)$ onto $H^{\frac12(\frac32)}(\Sigma _+)$. In other words, the domain of $L$ satisfies: 
$$
D(L)=H^{\frac12(\frac32)}(\Sigma _+)=\Lambda _+^{(-\frac12)}e^+\ol H^1(\Sigma _+^{\circ}),
\tag3.21
$$
which is contained in $ d^{\frac12}e^+\ol
H^1(\Sigma _+^\circ)+\dot H^{\frac32}(\Sigma _+)$.
\endproclaim 

\demo{Proof} It is seen from the second line in (3.19) that $S_3=L^{-1}-R$ is also
described by
$$
S_3=(RL_{1,+}+S_2)L^{-1}-R(L_{1,+}+L_{0,+})L^{-1}=S_2L^{-1}-RL_{0,+}L^{-1}.\tag3.22
$$
Here $S_2L^{-1}$ maps $\ol H^{\frac12}(\Sigma _+^\circ)$ into
$\E_\frac12 (\Sigma _+)$ in view of (3.19). For the
other term, we note that $L_{0,+}$ maps $\dot H^{1-\varepsilon
}(\Sigma _+)$ into $\ol H^{1-\varepsilon }(\Sigma _+^\circ)$, since
an extension by zero is understood, and $R$ maps the latter space into
$H^{\frac12(2-\varepsilon )}(\Sigma _+)$. Thus $S_3$ maps 
$\ol H^{\frac12}(\Sigma _+^\circ)$ into $H^{\frac12(2-\varepsilon
)}(\Sigma _+)$. Since $R$ maps 
$\ol H^{\frac12}(\Sigma _+^\circ)$ into $H^{\frac12(\frac32)}(\Sigma
_+)$, it follows that $L^{-1}$ maps $\ol H^{\frac12}(\Sigma _+^\circ)$ into
$H^{\frac12(\frac32)}(\Sigma_+)$. Thus $D(L)\subset
H^{\frac12(\frac32)}(\Sigma_+)$.

The opposite inclusion also holds,
since $r^+L_{1}$ maps $H^{\frac12(\frac32)}(\Sigma_+)$ into  $\ol
H^{\frac12}(\Sigma _+^\circ)$, and
$H^{\frac12(\frac32)}(\Sigma_+)\subset \dot H^\frac12(\comega)$ by
Lemma 2.4,
which $r^+L_0$ maps into  $\ol
H^{\frac12}(\Sigma _+^\circ)$.

This shows the identity.
The last statement follows from
(2.3). \qed

\enddemo
\example {Remark 3.5} 
By this information we can explain more precisely in which way $D(L)$, known
to be contained in $\dot H^{1-\varepsilon }(\Sigma _+)$,
reaches outside of $\dot H^1(\Sigma _+)$, namely by certain nontrivial
elements of  $d^\frac12 e^+\ol
H^{1}(\Sigma _+^{\circ})$ (lying in $ H^{\frac12(\frac32)}(\Sigma _+)$).

Consider the spaces
in
local coordinates, where $\Sigma $ and $\Sigma _+$ are replaced by ${\Bbb
R}^{n-1}$ and $ \crp^{n-1}$. As a typical element of $x_{n-1}^\frac12 e^+\ol
H^{1}(\rp^{n-1})$ lying in $H^{\frac12(\frac32)}(\crp^{n-1})$, we
can take 
$$
\varphi (x')=c\, x_{n-1}^\frac12 K_0\psi , \quad c  = \Gamma (\tfrac32)^{-1},\tag3.23
$$
 where $\psi (x'')\in H^\frac12({\Bbb R}^{n-2})$. Here $K_0$ is the
 Poisson operator from $H^{\frac12}({\Bbb R}^{n-2})$ to $\ol
 H^1(\rp^{n-1})$ solving
$$
(1-\Delta )\zeta (x')=0\text{ on }\rp^{n-1},\quad \gamma _0\zeta =\psi \text{ on }{\Bbb R}^{n-2},
$$
namely 
$$\zeta =K_0\psi =\F^{-1}_{\xi '\to x'}(\ang{\xi ''}+i\xi
_{n-1})^{-1}\hat\psi (\xi ''))=\F^{-1}_{\xi ''\to x''}(e^{-\ang{\xi ''}x_{n-1}}\hat\psi (\xi '')),
$$
and $\varphi (x')=c\, x_{n-1}^\frac12 \zeta (x')$.

To verify that $\varphi (x')\in H^{\frac12(\frac32)}(\crp^{n-1})$, we
recall from \cite{G13}, Sect.\ 5, that
 the special boundary operator $\gamma _{\frac12,0}\colon
H^{\frac12(\frac32)}(\crp^{n-1})\to H^\frac12({\Bbb R}^{n-2})$
defined there
satisfies
$$
\gamma _{\frac12,0}\varphi =c^{-1}\gamma _0(x_{n-1}^{-\frac12}\varphi
(x'))=\gamma _0\Xi _+^{\frac12}\varphi , \text{ with }
\Xi _+^{\mu}=\operatorname{OP}\bigl((\ang{\xi ''}+i\xi _{n-1})^{\mu}\bigr),   
$$
 and has the right inverse $K_{\frac12,0}$, where
$$
\varphi =K_{\frac12,0}\psi =\Xi _+^{-\frac12}e^+K_0\psi =c\, x_{n-1}^{\frac12}K_0\psi ,
$$
cf.\ \cite{G13}, Cor.\ 5.3, and the analysis in Th.\ 5.4 there. 

Now $\varphi $ defined by (3.23) is
not in $\dot H^1$ (nor in $\ol H^1$) near $x_{n-1}=0$, since 
$$
\partial_{x_{n-1}}\varphi (x')=\tfrac12 x_{n-1}^{-\frac12}\zeta (x') +x_{n-1}^{\frac12 }\partial _{x_{n-1}}\zeta (x'),
$$
where $x_{n-1}^{\frac12 }\partial _{x_{n-1}}\zeta (x')$ is clearly
$L_2$-integrable over ${\Bbb R}^{n-2}\times [0,1]$, but
$x_{n-1}^{-\frac12 }\zeta (x')$ is not so:
$$
\aligned
\int_{{\Bbb R}^{n-2}}&\int_{0<x_{n-1}<1}|x_{n-1}^{-\frac12}\zeta |^2\,
d{x_{n-1}}dx''
\\
&=(2\pi )^{2-n}\lim_{\delta \to 0}\int_{{\Bbb
R}^{n-2}}\int_{\delta <x_{n-1}<1}x_{n-1}^{-1}e^{-2\ang{\xi ''}x_{n-1}}|\hat \psi (\xi
'')|^2\, d{x_{n-1}}d\xi '' \\
&\ge (2\pi )^{2-n}\lim_{\delta \to 0}\int_{{\Bbb
R}^{n-2}}\int_{\delta <x_{n-1}<1}x_{n-1}^{-1}e^{-2\ang{\xi ''}}|\hat \psi (\xi
'')|^2\, d{x_{n-1}}d\xi '' \\
&= (2\pi )^{2-n}\lim_{\delta \to 0}|\log \delta |\int_{{\Bbb
R}^{n-2}}e^{-2\ang{\xi ''}}|\hat \psi (\xi
'')|^2\,d\xi '' 
=+\infty 
,\endaligned\tag3.24
$$
when $\psi \ne 0$. (It does not help to take $\psi $ very smooth.)
\endexample

We consequently have  for $D(A_{\chi ,\Sigma _+})$:

\proclaim{Corollary 3.6} The domain of $A_{\chi ,\Sigma _+}$ satisfies
$$
D(A_{\chi ,\Sigma _+})\subset D(A_\gamma )+K_\gamma
H^{\frac12(\frac32)}(\Sigma _+)\subset \ol H^2(\Omega )+K_\gamma
(e^+d(x')^{\frac12}\ol H^{1}(\Sigma _+^\circ))\tag3.25
$$
(where we recall that  $e^+$ denotes the extension from $\Sigma _+$ by zero on $\Sigma
_-$, and $d(x')$ is a $C^\infty $-function on $\Sigma _+$ proportional
to $\operatorname{dist}(x',\partial\Sigma _+)$ near $\partial\Sigma
_+$). 

All elements of $K_\gamma
H^{\frac12(\frac32)}(\Sigma _+)$ are reached from $D(A_{\chi ,\Sigma _+})$.

Nontrivial elements of  $K_\gamma
(e^+d(x')^{\frac12}\ol H^{1}(\Sigma _+^\circ))$ are reached, that are not
in $K_\gamma \dot H^1(\Sigma _+)$, nor in $K_\gamma (e^+\ol
H^1(\Sigma _+^\circ))$ (as in Remark {\rm 3.5}), hence not in $\ol
H^{\frac32}(\Omega )$.
.
\endproclaim

\demo{Proof} It is known from \cite{G68}, Th.\ II.1.2 that 
$$
D(A_{\chi ,\Sigma _+})\subset D(A_\gamma )\dot +D(T) = D(A_\gamma )\dot +K_\gamma
D(L),
$$
when we use that $A_\gamma =A_\beta $ and $K_\gamma D(L)=D(T)$ with
the notation used there. Here all elements of $D(T)$ are
reached, in the sense that for any $z\in D(T)$ there is a $v\in
D(A_\gamma )$ such that $u=v+z\in D(A_{\chi ,\Sigma _+})$. Since $D(L)=H^{\frac12(\frac32)}(\Sigma _+)$, this shows the
first inclusion in (3.25) and the first statement afterwards.

For the remaining part we use the last information in Theorem
3.4. Since  $K_\gamma \dot H^{\frac32}\subset \ol H^2(\Omega )$,
this implies the second inclusion in (3.25). Remark 3.5 shows how
nontrivial nonsmooth elements occur.\qed
\enddemo

\subhead 3.4 The spectrum of the Krein term  \endsubhead

The spectral asymptotic behavior of the Krein term
$$
M=A_{\chi ,\Sigma _=}^{-1}-A_\gamma ^{-1}=\inj_V\gamma
_V^{-1}L^{-1}(\gamma _V^{-1})^*\pr_V \tag 3.26
$$
will now be determined. We assume $n\ge 3$ in this section since  applications on $\Sigma $
of Laptev's result quoted in (2.14) requires the dimension $m$ to be $\ge
2$, i.e., $n-1\ge 2$. It is used to show that some cut-off terms have
a better asymptotic behavior than the one we are aiming for, hence can
be disregarded. (We believe that there are ways to handle the case
$n-1=1$, either by establishing weaker versions of (2.14), or by using the
variable-coefficient factorization of the principal symbol of $L$, but
we refrain from
making an effort here. The case $n=2$ was included in \cite{G11a} for
$A$ principally Laplacian.) 

First we study the spectrum of the factor $L^{-1}$.

\proclaim{Theorem 3.7} $S_3$ belongs to $\frak
S_{(n-1)/(\frac32 -\varepsilon ),\infty }$, and $L^{-1}$ belongs to $\frak
S_{n-1,\infty }$ (when the operators are extended to $L_2(\Sigma _+)$
by closure).

The eigenvalues of  $L^{-1}$ have the asymptotic behavior:
$$
\mu _j(L^{-1} )j^{1/(n-1)} \to c(L )^{1/(n-1)}\text{ for }j\to\infty ,\tag3.27
$$
where 
$$
c(L ) =
\tfrac1{(n-1)(2\pi )^{n-1}}\int_{\Sigma _+  }\int_{|\xi '
|=1}
\kappa _0(x',\xi ')^{-(n-1)}
\,d\omega (\xi ') dx '.\tag3.28
$$

\endproclaim

\demo{Proof} Recall that $L^{-1}$ acts as follows:  $$
L^{-1}=R+S_3=\Lambda _{+,+}^{(-\frac12)}(\widetilde Q_++G)\Lambda _{-,+}^{(-\frac12)}+S_3,\tag3.29
$$ 
cf.\ (3.18). By application of Theorem 2.6 to $R$ we find that the
singular values $s_j(R)$ behave as in (3.27)--(3.28), where the constant is
as in (3.28) since the principal pseudodifferential symbol of $R$ is
$\kappa _0^{-1}$. In particular,  $R\in \frak S_{n-1,\infty }$.
  
Since the closure of $L^{-1}$ maps $L_2(\Sigma _+)$
continuously into $\dot H^{\frac12 -\varepsilon }(\Sigma _+)$ by
Proposition 3.3, it belongs to
 $\frak S_{(n-1)/(\frac12-\varepsilon
),\infty }$. Moreover (cf.\ (3.19)), $S_1\in \bigcap_{\tau >0}\frak
S_{\tau ,\infty }$, and 
$L_{0,+}$ is bounded in $L_2(\Sigma _+)$. Then $L^{-1}S_1$ is in $\bigcap_{\tau >0}\frak S_{\tau ,\infty }$, and
$L^{-1}L_{0,+}R\in \frak S_{(n-1)/(\frac12 -\varepsilon ),\infty }\cdot
\frak S_{n-1,\infty }\subset \frak S_{(n-1)/(\frac32 -\varepsilon ),\infty
}$ by the rule (2.10), using that $S_1$ and $L_{0,+}R$ map into spaces
where $L^{-1}$ coincides with its $L_2$-closure. Therefore by (3.20),
$$
S_3\in \frak S_{(n-1)/(\frac32-\varepsilon ),\infty }.
$$
Now since $L^{-1}$ acts like $R+S_3$, its closure is in  $ \frak
S_{n-1,\infty }$.
This show the first statement in the theorem.

The last statement follows, since $S_3$ is of a better Schatten class
than $R$, so that (2.13) implies that the $L_2$-closure of $L^{-1}$ has the same asymptotic
behavior of singular values as $R$. Since $L^{-1}$ is symmetric in $L_2$, the
$L_2$-closure is selfadjoint, so its
singular values are eigenvalues; they are consistent with the
eigenvalues of $L^{-1}$ by
Proposition 3.3.
\qed
\enddemo

\comment
These properties of $L^{-1}$ were obtained considerably faster in \cite{G11a}, because we
now have formulations in terms of truly pseudodifferential operators
and can use their composition rules and spectral estimates.

Recalling the structure of $R$ we have that
$$
L^{-1}=
$$
where $S_3$ is of a better Schatten class than the first term, so it
will not contribute to the principal spectral asymptotics. 
\endcomment

We  now turn to the Krein term $M$ recalled in (3.26).
Proceeding as in \cite{G11a} Sect.\ 5.4, we have for the
eigenvalues:
$$
\mu _j(M)=\mu _j(\inj_V\gamma _V^{-1}L^{-1}(\gamma _V^{-1})^*\pr_V)
=\mu _j(L^{-1}(\gamma _V^{-1})^*\gamma _V^{-1})=\mu _j(L^{-1}P_{1,+}),
$$
where $P_1=K_\gamma ^*K_\gamma $ is a selfadjoint nonegative
invertible elliptic
$\psi $do of order $-1$; in view of (3.14) it has principal symbol
$(\kappa _++\bar\kappa _+)^{-1}=\ul a_{nn}(2\kappa _0)^{-1}$. Let
$P_2=P_1^{\frac12}$, then we continue the calculation as follows:
$$
\mu _j(M)=\mu _j(L^{-1}r^+P_2P_2e^+)=\mu _j(P_2e^+L^{-1}r^+P_2)=\mu _j(r^+P_2e^+L^{-1}r^+P_2e^++S_4),
$$
where $S_4$ is a sum of three terms, each one a product of $\psi $do's
and cutoff functions of a total order $-2$, and each containing a
factor either $r^-P_2e^+$ or $r^+P_2e^-$ (or both).  To the terms in
$S_4$ we can apply (2.14) together with 
product rules, concluding that they are
in $\frak S_{(n-1)/(2+\theta ),\infty }$ for some $\theta >0$. 

The operator (cf.\ (3.25))
$$
M_1=r^+P_2e^+L^{-1}r^+P_2e^+=P_{2,+}\Lambda _{+,+}^{(-\frac12)}(\widetilde Q_++G)\Lambda _{-,+}^{(-\frac12)}P_{2,+}+P_{2,+}S_3P_{2,+}
$$
is selfadjoint nonnegative, so its eigenvalues $\mu _j$ coincide with
the $s$-values. We can apply Theorem 2.6 to the first term, obtaining a spectral
asymptotic formula (2.16)--(2.17) with $t/n$ replaced by $2/(n-1)$;
then the addition of the second term which lies in a better weak
Schatten class  $\frak S_{(n-1)/(2+\theta ),\infty }$ preserves the formulas.

Finally $M$ (likewise selfadjoint nonnegative) differs from $M_1$ by the operator $S_4$ in a better weak
Schatten class, so the spectral asymptotic formula carries over to
this operator.

Hereby we obtain the theorem:

\proclaim{Theorem 3.8} The eigenvalues of $M=A_{\chi ,\Sigma
_+}^{-1}-A_\gamma ^{-1}$ have the asymptotic behavior:
$$
\mu _j(M )j^{2/(n-1)} \to c(M )^{2/(n-1)}\text{ for }j\to\infty ,\tag3.30
$$
where 
$$
c(M ) =
\tfrac1{(n-1)(2\pi )^{n-1}}\int_{\Sigma _+  }\int_{|\xi '
|=1}
\Big(\frac{\ul a_{nn}(x')}{2\kappa _0(x',\xi ')^2}
\Big)^{(n-1)/2}
\,d\omega (\xi ') dx '.
\tag3.31
$$

\endproclaim

\demo{Proof} It remains to account for the value of the constant
$c(M)$. It follows, since $P_2^2=P_1$ has principal symbol $\ul
a_{nn}(2\kappa _0)^{-1}$ and the $\psi $do part of $L^{-1}$ has
principal symbol $\kappa _0^{-1}$. \qed
\enddemo

\example{Remark 3.9}
We take the opportunity to recall two corrections to \cite{G11a}
(already mentioned in \cite{G11b}): Page 351, line 4 from below, delete
``$H^{\frac12 }(\Sigma _+^\circ)\subset$'', replace ``$H^1(\Sigma )$'' by
``$L_2(\Sigma )$''. Page 361, line 4, replace ``(Th.\ 3.3)'' by ``(Th.\ 4.3)''. 
\endexample

\enddocument

\end

%% file: gmacro2.tex
\def\supp{\operatorname{supp}}

\def\pr{\operatorname{pr}}

\def\pr{\operatorname{pr}}

\def\crp{\overline{\Bbb R}_+}

\def\rnp{{\Bbb R}^n_+}
\def\rnm{\Bbb R^n_-}
\def\rnpm{\Bbb R^n_\pm}
\def\crnp{\overline{\Bbb R}^n_+}

\def\crnpm{\overline{\Bbb R}^n_\pm}
\def\comega{\overline\Omega }

\def\ang#1{\langle {#1} \rangle}

\def\Pfrac{\tsize\frac1{\raise 1pt\hbox{$\scriptstyle p$}}}
\def\pfrac{\frac1{\raise 1pt\hbox{$\scriptscriptstyle p$}}}
\def\Pfracc#1{\tsize\frac{#1}{\raise 1pt\hbox{$\scriptstyle p$}}}
\def\pfracc#1{\frac{#1}{\raise 1pt\hbox{$\scriptscriptstyle p$}}}

\def\simto{\overset\sim\to\rightarrow}

\def\Zfrac{\tsize\frac1{\raise 1pt\hbox{$\scriptstyle z$}}}
\def\zfrac{\frac1{\raise 1pt\hbox{$\scriptscriptstyle z$}}}

\def\rp{ \Bbb R_+}

\def\Ama{A_{\max}}
\def\inj{\text{\rm i}}

\def\R{\Bbb R}

\def\ol{\overline}
\def\SD{\Cal S}
\def\E{\Cal E}
\def\F{\Cal F}
\def\D{\Cal D}